\documentclass[11pt]{amsart}
\usepackage{mathrsfs}
\usepackage{amsfonts}
\usepackage{amsthm}
\usepackage{amssymb}
\usepackage{latexsym}
\usepackage{color}
\usepackage{xcolor}
\usepackage{soul}
\newtheorem{theorem}{Theorem}
\newtheorem{lemma}{Lemma}
\newtheorem{corollary}{Corollary}

\begin{document}

\title[regularity in the $\overline{\partial}$--Neumann problem]{Regularity in the $\overline{\partial}$--Neumann problem, D'Angelo forms, and Diederich--Forn\ae ss index}
\author{Emil J. Straube}
\thanks{Work supported in part by NSF grant DMS--2247175. This article is an expanded version of a lecture first given at the Schr\"{o}dinger Institute in Vienna during the workshop on Analysis and Geometry in Several Complex Variables, Nov. 20--24, 2023.}

\address{Department of Mathematics, Texas A\&M University, College Station, Texas, USA}
\email{e-straube@tamu.edu}
\date{April 3, 2025}

\maketitle

\begin{center}
\emph{Dedicated to the memory of Joseph~J.~Kohn}
\end{center}

\begin{abstract}
This article chronicles a development that started around 1990 with \cite{BoasStraube91}, where the authors showed that if a smooth bounded pseudoconvex domain $\Omega$ in $\mathbb{C}^{n}$ admits a defining function that is plurisubharmonic at points of the boundary, then the $\overline{\partial}$--Neumann operators on $\Omega$ preserve the Sobolev spaces $W^{s}_{(0,q)}(\Omega)$, $s\geq 0$. The same authors then proved a further regularity result and made explicit the role of D'Angelo forms for regularity (\cite{BoasStraube93}). A few years later, Kohn (\cite{Kohn99}) initiated a quantitative study of the results in \cite{BoasStraube91} by relating the  Sobolev level up to which regularity holds to the Diederich--Forn\ae ss index of the domain. Many of these ideas were synthesized and developed further by Harrington (\cite{Harrington11,Harrington19,Harrington22}). Then, around 2020, Liu (\cite{Liu19b, Liu19}) and Yum (\cite{Yum21}) discovered that the DF--index is closely related to certain differential inequalities involving D'Angelo forms. This relationship in turn led to a recent new result which supports the conjecture that DF--index one should imply global regularity in the $\overline{\partial}$--Neumann problem (\cite{LiuStraube22}). Much of the work described above relies heavily on Kohn's groundbreaking contributions to the regularity theory of the $\overline{\partial}$--Neumann problem.
\end{abstract}

\section{A brief history}\label{history} 

 The $L^{2}$--Sobolev theory of the $\overline{\partial}$--Neumann problem as we know it would be unthinkable without Kohn's many foundational contributions. In this short history, I will concentrate on the theory on domains in $\mathbb{C}^{n}$ with scalar valued forms, as it relates to the items in the title. For the theory on domains in complex manifolds and bundle valued forms, the reader should consult the excellent surveys \cite{Berndtsson10, Varolin10, McNealVarolin15} and the book \cite{Demailly12}, and their references.\footnote{This brief history is essentially an abridged version of chapter 1 in \cite{Straube10}, updated to account for developments after 2010 related to the themes in the title.}

The $\overline{\partial}$--Neumann problem is the problem of inverting the complex Laplacian $\Box = \overline{\partial}\overline{\partial}^{*}+\overline{\partial}^{*}\overline{\partial}$, say on a bounded pseudoconvex domain in $\mathbb{C}^{n}$. Here the adjoint is with respect to the $L^{2}$--norm of forms. Being in the domain of the adjoint $\overline{\partial}^{*}$ involves a boundary condition (as usual, from integration by parts), namely the normal component of a form $u$ should vanish on the boundary. Inverting $\Box=\overline{\partial}\overline{\partial}^{*}+\overline{\partial}^{*}\overline{\partial}$ thus comes with two boundary conditions, one for $u$ itself, and one for $\overline{\partial}u$. The latter one is analogous to the boundary condition for the Neumann problem for the Laplacian: the normal derivative of the function $f$, i.e. the normal component of the differential $df$, should vanish on the boundary. While the complex $\overline{\partial}\oplus\overline{\partial}^{*}$ is elliptic, so that there is no problem with interior regularity of $\Box$, the $\overline{\partial}$--Neumann boundary conditions turn out to be noncoercive, so that the classical theory of elliptic boundary value problems does not apply. This made the problem difficult.

The basic estimate for the $\overline{\partial}$--Neumann problem is the $L^{2}$--estimate \eqref{basic} below. An equivalent statement is that $\Box$ has a bounded inverse, $N$, in $L^{2}$. These estimates were proved in the sixties by Kohn for strictly pseudoconvex domains (\cite{Kohn61, Kohn63, Kohn63a, Kohn64}) (see also Morrey \cite{Morrey63, Morrey64}), and by H\"{o}rmander (\cite{Hormander65, Hormander73}) in general (see also \cite{AV65}). The key to these results, and thus to the $L^{2}$--theory of the $\overline{\partial}$--Neumann problem, is the so called Kohn--Morrey--H\"{o}rmander formula (\cite{Kohn63, Kohn64, Morrey64, Hormander65}) and a generalization due to Ohsawa--Takegoshi (\cite{Ohsawa88, OhsawaTakegoshi87}), McNeal (\cite{McNeal96, McNeal05}), and Siu (\cite{Siu96}). Kohn actually showed that in the strictly pseudoconvex case, $N$ gains one derivative in the Sobolev scale; because the gain is less than the order of the operator, one says the operator is subelliptic. This raised the question under what circumstances $N$ is subelliptic. After Kohn (\cite{Kohn79}) gave sufficient conditions and noted that work of Diederich and Forn\ae ss (\cite{DF78}) implies that these conditions are satisfied on domains with real analytic boundaries, the question was settled by his students Catlin and D'Angelo in the eighties (\cite{Catlin83, Catlin84a, Catlin87b, D'Angelo87, DA80, DA82}): $N$ is subelliptic if and only if the domain is of finite (D'Angelo) type. When $N$ is subelliptic, it is in particular compact as an operator on $L^{2}$, but compactness is a much weaker property. Earlier work of Kohn and Nirenberg (\cite{KohnNirenberg65}) implies that then, $N$ is continuous (in fact also compact) in Sobolev norms. They also developed the technique of elliptic regularization to deal with the problem of \emph{a priori} estimates. I do not discuss subellipticity or compactness of $N$ here, but instead refer the reader to the surveys (\cite{D'AngeloKohn99, D'AngeloKohn10, D'Angelo10, FuStraube99}). For efforts to combine the approach to global regularity via compactness with the methods to be discussed below, see \cite{Straube05, Harrington11}. 

When $N$ is subelliptic or compact, the relevant commutators with vector fields can be absorbed not by virtue of properties of the vector fields one uses to control derivatives, but because estimates better than \eqref{basic} hold. When only \eqref{basic} is available, something good has to come from the vector fields. For example, if $X$ is a $(1,0)$--field with holomorphic coefficients, then $[\overline{\partial}, X] = 0$. With this observation as a starting point, there appeared throughout the eighties and into the early nineties a number of results that concerned domains with `transverse symmetries' (\cite{BellBoas81, Barrett82, Barrett86, Straube86, Chen89}) or variants (\cite{BCS88, BoasStraube89, BoasStraube92, Chen87}), in particular Reinhardt domains and many circular domains. Many of these regularity results were actually proved for the Bergman projection, but in \cite{BoasStraube90}, Boas and I proved that Sobolev estimates for the $\overline{\partial}$--Neumann operator $N_{q}$ are equivalent to Sobolev estimates for the three Bergman projections at consecutive levels, $P_{q-1}$, $P_{q}$, and $P_{q+1}$, $1\leq q\leq n$. This result relies heavily on another seminal contribution by Kohn, namely his weighted theory (\cite{Kohn73}), as well as on what is now known as Kohn's formula $P_{q}=id -\overline{\partial}^{*}N_{q+1}\overline{\partial}$ that expresses the Bergman projection in terms of the $\overline{\partial}$--Neumann operator (\cite{Kohn63, Kohn64}, p. 235 in \cite{Kohn77}, (14) in \cite{Kohn84}). Around 1990, these techniques, essentially relying on vanishing commutators, had gone as far as they could, and one had to deal with actual (i.e. nonvanishing) commutators. In \cite{BoasStraube91}, three observations were made: (1) derivatives of type $(0,1)$ and tangential derivatives of type $(1,0)$ are benign for the $\overline{\partial}$--Neumann problem, (2) commutators with strictly pseudoconvex directions and commutators with the complex normal can always be adjusted for free\footnote{Though not explicitly stated in \cite{BoasStraube91} in this way, it is implicit, see also \cite{BoasStraube93}. This idea had antecedents going back at least to \cite{DerridjTartakoff76, Komatsu76}, in the context of nondegenerate Levi form and real analytic boundary.}, and, crucially, (3) it suffices to have $(1,0)$--vector fields that are (close to the) complex normal and with favorable `\emph{approximate} commutator conditions'. The notion of favorable `approximate commutator condition' is implemented via a one parameter \emph{family} of vector fields with better and better commutation properties, see Theorem \ref{vectorfield}. It allows for considerable flexibility, while the corresponding exact condition, with only one vector field, is rather restrictive.\footnote{See \cite{Derridj91}, Th\'{e}or\`{e}me 2.6 and the associated remark, and Remark 2 in section 4 of \cite{BoasStraube93}.} In particular, we showed in \cite{BoasStraube91} that these conditions are satisfied when the domain admits a defining function that is plurisubharmonic at points of the boundary.\footnote{We originally circulated a preprint on convex domains only that was inspired by Alan Noell's observation that the normal is constant along a Levi null curve in the boundary of a convex domain in $\mathbb{C}^{n}$ (\cite{Noell91}). So-Chin Chen obtained regularity for convex domains in $\mathbb{C}^{2}$ independently at about the same time in \cite{Chen89b}. A little earlier, Bonami and Charpentier (\cite{BonamiCharpentier88}) had shown that a plurisubharmonic defining function implies regularity in $W^{1/2}$.} In light of observation (1) above, the commutator condition is on the normal $(1,0)$--component of the commutator of $\overline{L}$ with the complex normal, where $L$ is a tangential $(1,0)$--field (standing in for $\overline{\partial}$). A simple computation expresses this component in terms of the complex Hessian of a defining function, whence the importance of plurisubharmonicity of the defining function. The proof in \cite{BoasStraube91} that (3) above suffices to obtain Sobolev estimates relies again in essential ways on Kohn's work. An argument quite different form that in \cite{BoasStraube91} allowed Gallagher and McNeal (\cite{Herbig-McNeal05}) in 2006 to extended the main result in \cite{BoasStraube91} so that the assumption on the defining function is tailored to the form level.

Essentially concurrently\footnote{In fact, Barrett and I gave back to back talks on our results at a workshop at Oberwolfach in early summer of 1990.} with \cite{BoasStraube91}, Barrett showed that, contrary to many people's expectations at the time, Sobolev estimates for the $\overline{\partial}$--Neumann operator do not always hold: they fail on the Diederich--Forn\ae ss worm domains (\cite{Barrett92}; see \cite{Barrett84} and \cite{Kisel91} for antecedents) when the Sobolev level is at least $\pi/(\text{total amount of winding})$.  A few years later, Christ strengthened Barrett's result by showing that the $\overline{\partial}$--Neumann operator on any worm domain $\Omega$ does not even preserve $C^{\infty}_{(0,1)}(\overline{\Omega})$ (\cite{Christ96}; compare also \cite{Siu00}).

Enter a family of one--forms introduced into the literature by D'Angelo in \cite{DA79, D'Angelo87}, who had observed that the relevant commutator component can also be expressed as a one--form given by a Lie derivative naturally associated to the (boundary of the) domain. Since there are choices involved, there actually is a family of forms. However, the relevant properties are independent of the particular form in the family. The role of these forms was made explicit in \cite{BoasStraube93}, see also \cite{Straube10}, section 5.9. Taking advantage of the simple expression as a Lie derivative, so that formulas form the calculus of differential forms can be exploited, the authors in \cite{BoasStraube93} discovered an important property of these forms: their differential is closed on the null space of the Levi from (\cite{BoasStraube93}, \cite{Straube10}, Lemma 5.14). In particular, the restriction of a D'Angelo form to a complex manifold in the boundary defines a DeRham cohomology class there. When this class vanishes, the form is exact on this manifold and thus, if the boundary is otherwise of finite type, exact on the null space of the Levi form. This latter condition is sufficient for the existence a family of `good' vector fields as in the previous paragraph (\cite{BoasStraube91, BoasStraube93}), and global (exact) regularity of the $\overline{\partial}$--Neumann operator follows. The appearance of this cohomology class explains in particular why an analytic disc in the boundary is always benign for global regularity (\cite{BoasStraube92}), while an analytic annulus is bad for the worm domains (\cite{Barrett92, Christ96}), but is benign for `nowhere wormlike' Hartogs domains in $\mathbb{C}^{2}$(\cite{BoasStraube92}). When there is a foliation by complex submanifolds in the boundary, say of co--dimension one, whether the form is exact on the leaves of the foliation by a `global' function is equivalent to a question much studied in foliation theory, namely whether the foliation can be defined globally by a closed one--form (\cite{Straube06}, Section 3.1, \cite{Straube10}, Lemma 5.18). Regularity up to a finite level depending on how close to exact the D'Angelo forms are was studied by Harrington and Liu (\cite{HarringtonLiu20}).

In \cite{Kohn99}, Kohn pioneered a quantitative study of the regularity estimates in \cite{BoasStraube91}, in terms of `how close to having a plurisubharmonic defining function' the domain is, in other words, in terms of the Diederich--Forn\ae ss index (see also \cite{PinZamp14}). He derived Sobolev regularity up to a level depending on the index. Assuming some control on $\nabla h_{\eta}$, where $\rho_{\eta}=e^{h_{\eta}}\rho$ are the defining functions in the definition of the index, he obtained regularity at all levels if the index is one\footnote{Berndtsson and Charpentier (\cite{BernChar00}) and McNeal (\cite{McNealUN}) obtained regularity for levels less than half of the DF--index in general, see also Michel and Shaw (\cite{MichelShaw01}) for index one; an important point in these papers is that they only require the domain to have Lipschitz boundary.}. Many of the above ideas were synthesized and developed further by Harrington (\cite{Harrington11, Harrington19, Harrington22}). These works seemed to indicate that index one should imply regularity at all levels\footnote{Note that index one does not imply the existence of a plurisubharmonic defining function, not even locally and just at boundary points (\cite{Behrens85, GallagherHarz22}), so \cite{BoasStraube91} does not apply.}. An important step in this direction was achieved by Bingyuan Liu (\cite{Liu19, Liu19b}). He expressed the DF--index as an index on the boundary and derived certain differential inequalities with which to characterize the index (see \eqref{DFDA3} below). This allowed him to prove the first result relating index one to regularity that does not require any \emph{a priori} control on $h_{\eta}$ or $\nabla h_{\eta}$ ($h_{\eta}$ as above); instead, he assumes a comparability condition on the Levi eigenvalues and, in addition, a somewhat restrictive geometric condition on the set of infinite type points. Under these conditions, he was able to show that index one implies (exact) regularity. Using his characterization of the index, Liu was also able to compute the index of the Diederich--Forn\ae ss worm domains (\cite{Liu19}); it equals $\pi/(\text{total winding}+\pi)$. Jihun Yum (\cite{Yum21}) then realized that Liu's characterization of the index can be reformulated in terms of certain estimates satisfied by D'Angelo's forms and their differentials (see Theorem \ref{DFDA2} below). Dall'Ara and Mongodi (\cite{Dall'AraMongodi21}) subsequently reworked and clarified these ideas. The estimates evoke the $L^{2}$--type sufficient condition for regularity introduced earlier in \cite{Straube05}. This connection led to a proof with Liu (\cite{LiuStraube22}) that one can dispense with the ad hoc geometric condition in \cite{Liu22}: when the Levi eigenvalues are comparable, then trivial DF--index implies Sobolev estimates of all orders for the $\overline{\partial}$--Neumann operator. In particular, for domains in $\mathbb{C}^{2}$, trivial index implies regularity of all orders. I expect this conclusion  to hold in all dimensions, without the comparability assumption on the Levi eigenvalues, but for now, that remains open. 

\medskip

In the remainder of the paper, I will discuss in more detail the interactions between D'Angelo forms, regularity, and the DF--index. Section 2  gives the needed definitions and notation along with some basic results. Section \ref{alphaest} shows how properties of D'Angelo forms impact regularity, section \ref{DA-DF} expresses the index in terms of a differential estimate for D'Angelo forms, and section \ref{DF-reg} discusses trivial DF--index and regularity.

\section{Definitions}\label{def}

Unless other wise stated, $\Omega$ will be a bounded smooth ($C^{\infty}$) pseudoconvex domain in $\mathbb{C}^{n}$ throughout.

\smallskip
\underline{\emph{2.1. The $\overline{\partial}$--Neumann problem:}}
Proofs and more details for this subsection can be found in \cite{ChenShaw01, Straube10}.           
The $L^{2}$--space of $(0,q)$--forms, denoted by $L^{2}_{(0,q)}(\Omega)$, consists of forms 
\begin{equation}\label{def}
u=\sideset{}{'}\sum_{|J|=q}u_{J}d\overline{z_{J}},\;\text{with}\;u_{J}\in L^{2}(\Omega)\;,  d\overline{z_{J}}=d\overline{z_{j_{1}}}\wedge \cdots\wedge d\overline{z_{j_{q}}}\;, 
\end{equation}
and the prime denoting summation over strictly increasing $q$--tuples. Here, $0\leq q\leq n$. The inner product is
\begin{equation}\label{inner}
 (u, v) = \big(\sideset{}{'}\sum_{|J|=q}u_{J}d\overline{z_{J}}, \sideset{}{'}\sum_{|J|=q}v_{J}d\overline{z_{J}}\big) = \sideset{}{'}\sum_{|J|=q}\int_{\Omega}u_{J}(z)\overline{v_{J}(z)}\,dV(z)\;.
\end{equation}
We have
\begin{equation}\label{dbar}
 \overline{\partial}u = \sum_{j=1}^{n}\sideset{}{'}\sum_{|J|=q}\frac{\partial u_{J}}{\partial\overline{z_{j}}}d\overline{z_{j}}\wedge d\overline{z_{J}}\;;
 \end{equation}
the domain of $\overline{\partial}$ consists of all forms such that when $\overline{\partial}u$ is computed in the sense of distributions, the resulting form is in $L^{2}_{(0,q+1)}(\Omega)$. $\overline{\partial}$ is closed and densely defined, and so has a Hilbert space adjoint $\overline{\partial}^{*}$. For $u\in dom(\overline{\partial}^{*})\subset L^{2}_{(0,q+1)}(\Omega)$, this adjoint is
\begin{equation}\label{dbar*}
\overline{\partial}^{*}u = -\sideset{}{'}\sum_{|K|=q-1}\big(\sum_{j=1}^{n}\frac{\partial u_{jK}}{\partial z_{j}}\big)d\overline{z_{K}}\;.
 \end{equation}
The complex Laplacian is defined by
\begin{equation}\label{laplace}
 \Box_{q} u = (\overline{\partial}\overline{\partial}^{*} + \overline{\partial}^{*}\overline{\partial})u\;,
 \end{equation}
with domain so that the two compositions are defined. $\Box_{q}$ is densely defined and self--adjoint. The $\overline{\partial}$--Neumann problem consists of inverting $\Box_{q}$, with estimates. 

 The crucial fact about $\overline{\partial}$ on bounded smooth pseudoconvex domains is that it has closed range. Hence so do $\overline{\partial}^{*}$ and $\Box_{q}$. Specifically, there is the estimate
\begin{equation}\label{basic}
  \|u\|^{2} \leq \frac{D^{2}e}{q}\big(\|\overline{\partial}u\|^{2}+\|\overline{\partial}^{*}u\|^{2}\big)\;; u\in dom(\overline{\partial})\cap dom(\overline{\partial}^{*})\;. 
\end{equation}
                                                                                                                                                                                                                                Here $D$ is the diameter of $\Omega$, and $q$ is the form level. The space $dom(\overline{\partial})\cap dom(\overline{\partial}^{*})$ provided with the graph norm is a Hilbert space. Then \eqref{basic} says that the embedding $j_{q}: dom(\overline{\partial})\cap dom(\overline{\partial}^{*}) \rightarrow L^{2}_{(0,q)}(\Omega)$ is continuous. Because $\Box$ has trivial kernel (from \eqref{laplace}, \eqref{basic}) and is self--adjoint, it has dense range. Because this range is closed, it is all of $L^{2}_{(0,q)}(\Omega)$. It follows that $\Box_{q}$ has a bounded inverse. This inverse is the $\overline{\partial}$--Neumann operator $N_{q}$. It is not hard to see that $N_{q}=j_{q}\circ j_{q}^{*}$. $N_{q}$ commutes with $\overline{\partial}$ and $\overline{\partial}^{*}$ (for forms in the appropriate domains). If $u$ is a $\overline{\partial}$--closed $(0,q)$--form, this gives
 \begin{equation}\label{canonical}
 u = \overline{\partial}\overline{\partial}^{*}N_{q}u + \overline{\partial}^{*}\overline{\partial}N_{q}u =   \overline{\partial}(\overline{\partial}^{*}N_{q})u \;;
\end{equation}
 that is, $v=\overline{\partial}^{*}N_{q}u$ provides a solution to the equation $\overline{\partial}v=u$.  It is the (unique) solution with minimal $L^{2}$--norm and is referred to as the Kohn, or canonical, solution.

                                                                                                                                                        For $s\in\mathbb{R}$, we denote by $W^{s}(\Omega)$ the usual $L^{2}$--Sobolev spaces of order $s$, with norm $\|\cdot\|_{s}$, and by $W^{s}_{(0,q)}(\Omega)$ the corresponding spaces of $q$--forms (with the natural inner product induced by $W^{s}(\Omega)$). We say that $N_{q}$ is regular in Sobolev norms if there are estimates 
\begin{equation}\label{regular}
 \|N_{q}u\|_{s} \leq C_{s}\|u\|_{s}\;;\;s\geq 0\;.
\end{equation}
When one has Sobolev estimates as in \eqref{regular} for $N_{q}$, the same estimates hold for the Kohn solution operator $\overline{\partial}^{*}N_{q}$. That is, one can solve the $\overline{\partial}$--equation with estimates. 

A useful fact is that roughly speaking, barred derivatives and tangential unbarred derivatives are controlled by $\overline{\partial}u$ and $\overline{\partial}^{*}u$:
\begin{equation}\label{benign1}
 \sum_{j,J}\left\|\frac{\partial u_{J}}{\partial\overline{z_{j}}}\right\|_{k-1}^{2} \leq C\left(\|\overline{\partial}u\|_{k-1}^{2}+\|\overline{\partial}^{*}u\|_{k-1}^{2}+\|u\|_{k-1}^{2}\right)
\end{equation}
and
\begin{equation}\label{benign2}
  \|Yu\|_{k-1}^{2} \leq C\left(\|\overline{\partial}u\|_{k-1}^{2}+\|\overline{\partial}^{*}u\|_{k-1}^{2}+\|u\|_{k-1}\,\|u\|_{k}\right) \;,       
\end{equation}
where $Y\in T^{1,0}(b\Omega)$ is smooth on $\overline{\Omega}$. For a proof, see \cite{BoasStraube91}, estimates (2) and (3); \cite{Straube10}, Lemma 5.6.

 \smallskip

The observations about commutators of vector fields with strictly pseudoconvex directions and approximately favorable commutator conditions mentioned in section \ref{history} are contained in the following theorem from \cite{BoasStraube93} (see also \cite{Straube10}, Theorem 5.9), which is the starting point for the story in this survey. Denote by $K$ the set of boundary points of infinite type.
\begin{theorem}\label{vectorfield}
Let $\Omega$ be a smooth bounded pseudoconvex domain in $\mathbb{C}^{n}$ with smooth defining function $\rho$. Suppose there is a positive constant $C$ such that the following holds. For every positive $\varepsilon>0$ there is a vector field $X_{\varepsilon}$ of type $(1,0)$ whose coefficients are smooth in a neighborhood $U_{\varepsilon}$ of $K$, such that

(i) $|\arg X_{\varepsilon}\rho| < \varepsilon$ and $C^{-1} < |X_{\varepsilon}\rho| < C$ on $K$; and

(ii) $|\partial\rho([X_{\varepsilon}, \overline{L}])(P)| < \varepsilon$ whenever $L$ is a tangential $(1,0)$--field of unit length near $P\in K$ with $L(P)\in\mathcal{N}_{P}$.

\noindent Then the $\overline{\partial}$--Neumann operators $N_{q}$, $1\leq q\leq n$, are continuous in $W_{(0,q)}^{s}(\Omega)$ when $s\geq 0$.
\end{theorem}
\noindent Note that whether such a family of vector fields exists does not depend on the defining function $\rho$: near $b\Omega$, any other defining function is of the form $e^{h}\rho$. 
The reason why conditions are only need at (near) points of $K$ is that near points of finite type (very strong) subelliptic pseudolocal estimates hold (\cite{Catlin87b}) and give the necessary control. There is an interesting variant of Theorem \ref{vectorfield}
in \cite{HarringtonLiu20}: when $X_{\varepsilon}\rho$ is assumed real, then the uniform bounds are not needed, at the expense of having to consider commutators with all $L\in T_{P}^{1,0}(b\Omega)$, instead of just $L\in\mathcal{N}_{P}$, see Theorem \ref{HL20} below.

\smallskip

To conclude this subsection, consider a regularity property that on its face is weaker than continuity in $W^{s}_{(0,q)}(\Omega)$ for $s\geq 0$ (exact regularity), namely continuity in $C^{\infty}_{(0,q)}(\overline{\Omega})$ (global regularity). It is rather intriguing that no examples are known where the $\overline{\partial}$--Neumann operators are globally regular without being exactly regular.

\medskip

\underline{\emph{2.2. D'Angelo forms:}}
Next are the forms now commonly known as D'Angelo forms, introduced into the literature in \cite{DA79, D'Angelo87}. Let $\rho$ be a defining function for $\Omega$, set $(L_{n})_{\rho}=(1/|\partial\rho|^{2})\sum_{j=1}^{n}(\partial\rho/\partial\overline{z_{j}})\partial/\partial z_{j}$ (so that $(L_{n})_{\rho}\rho\equiv 1$), $T_{\rho}=(L_{n})_{\rho} - \overline{(L_{n})_{\rho}}$, and $\sigma_{\rho}=(1/2)(\partial\rho - \overline{\partial}\rho)$ (so that $\sigma_{\rho}(T_{\rho})\equiv 1$)\footnote{All quantities are well defined near $b\Omega$.}. Then we define
\begin{equation}\label{alpha1}
\alpha_{\rho} = -\mathcal{L}_{T_{\rho}}\sigma_{\rho}\;,
\end{equation}
where $\mathcal{L}$ denotes the Lie derivative; note that $\alpha$ is real. When the $(1,0)$--field $L$ belongs to $T^{1,0}(b\Omega)$, the definition of the Lie derivative and the fact that $\sigma_{\rho}$ annihilates $T^{1,0}(b\Omega)\oplus T^{0,1}(b\Omega)$ give
\begin{multline}\label{alpha2}
 \alpha_{\rho}(\overline{L}) = -T_{\rho}(\sigma_{\rho}(\overline{L}))+\sigma_{\rho}([T_{\rho}, \overline{L}]) = \sigma_{\rho}([T_{\rho}, \overline{L}]) \\
 = \frac{1}{2}(\partial\rho-\overline{\partial}\rho)([(T_{\rho}, \overline{L}])=\partial\rho([T_{\rho}, \overline{L}])=\partial\rho([(L_{n})_{\rho}, \overline{L}])
\end{multline}
($[T_{\rho}, \overline{L}]$ is tangential, so is annihilated by $(\partial\rho+\overline{\partial}\rho)$). Thus $\alpha_{\rho}$ picks out the normal $(1,0)$--components of the commutators $[T_{\rho}, \overline{L}]$ and $[(L_{n})_{\rho}, \overline{L}]$. We also have 
\begin{equation}\label{alpha3}
\alpha_{\rho}(L) = \overline{\alpha_{\rho}(\overline{L})} =\overline{\partial}\rho([-T_{\rho}, L]) = -\partial\rho([-T_{\rho}, L]) = \partial\rho([T_{\rho}, L])\;,
\end{equation}
again because $[-T_{\rho}, L]$ is tangential. For a \emph{constant} multiple of $(L_{n})_{\rho}$, \eqref{alpha2} and \eqref{alpha3} also hold, by a similar computation, using that $\sigma_{\rho}\big(\overline{(L_{n})_{\rho}}\big)$ is constant, and that $\big[T_{\rho}, \overline{(L_{n})_{\rho}}\big]$ is tangential:
\begin{equation}\label{L_ncomm}
 \alpha_{\rho}\big(\overline{(L_{n})_{\rho}}\big) = \partial\rho\big(\big[T_{\rho},\overline{(L_{n})_{\rho}}\big]\big)\;\;;\;\;  \alpha_{\rho}\big((L_{n})_{\rho}\big) = \partial\rho\big(\big[T_{\rho},(L_{n})_{\rho}\big]\big)\;.
\end{equation}

Alternatively, the commutator components in \eqref{alpha2} and \eqref{alpha3} can also be expressed in terms of the Hessian of the defining function.  For example, if $L=\sum_{k=1}^{n}w_{k}(\partial/\partial z_{k})$, then
\begin{multline}\label{alternative}
 \alpha_{\rho}(\overline{L})=\partial\rho([(L_{n})_{\rho}, \overline{L}])=-\sum_{j=1}^{n}\overline{L}\left(\left((L_{n})_{\rho}\right)_{j}\right)\frac{\partial\rho}{\partial z_{j}} \\
 = \sum_{j=1}^{n}\left((L_{n})_{\rho}\right)_{j}\overline{L}\left(\frac{\partial\rho}{\partial z_{j}}\right)= \sum_{j,k=1}^{n}\frac{\partial^{2}\rho}{\partial z_{j}\partial\overline{z_{k}}}((L_{n})_{\rho})_{j}\overline{w_{k}}\;.\;\;\;\;\;
\end{multline}
For the third equality, note that $\sum_{j=1}^{n}\left((L_{n})_{\rho})\right)_{j}(\partial\rho/\partial z_{j})$ is constant.

The commutators in \eqref{alpha2} -- \eqref{L_ncomm} are the relevant commutators for the $\overline{\partial}$--Neumann problem. First, because barred derivatives and complex tangential derivatives are under control by \eqref{benign1} and \eqref{benign2}, one only needs to consider $T_{\rho}$ derivatives. Second, one then only needs to estimate the normal $(1,0)$--component of the commutators $[\overline{\partial},(T_{\rho})^{k}]$ and $[\overline{\partial}^{*},(T_{\rho})^{k}]$, for the same reason. Suppose a form $v=\sideset{}{'}\sum_{|J|=q}v_{J}\overline{\omega_{J}}$ is supported in a special boundary chart, with the usual orthonormal basis $L_{1}, \cdots, L_{n}$ and dual basis $\omega_{1}, \cdots, \omega_{n}$, and differential operators acting as usual coefficientwise in this chart. Then these  considerations together with a computation, using \eqref{alpha2} -- \eqref{L_ncomm} for the normal components of the commutators, give the case $k=1$ of the following formulas (\cite{HarringtonLiu20}, Lemmas 3.1, 3.2; \cite{LiuStraube22}, Lemmas 3, 4): 
\begin{equation}\label{dbarcomm}
\big[\overline{\partial},(T_{\rho})^{k}\big]v
=-k\sideset{}{'}\sum_{|J|=q}\sum_{j}\alpha_{\rho}(\overline{L_{j}})\big((T_{\rho})^{k}v_{J}\big)(\overline{\omega_{j}}\wedge\overline{\omega_{J}})+A(v)+B(v)\;,\;\;\;\;\;
\end{equation}
and, for $v\in dom(\overline{\partial}^{*})$,
\begin{equation}\label{dbar*comm}
 \big[\overline{\partial}^{*},(T_{\rho})^{k}\big]v = -k\sideset{}{'}\sum_{|S|=q-1}\bigr(\sum_{j}\alpha_{\rho}(L_{j})\big((T_{\rho})^{k}v_{jS}\big)\bigr)\overline{\omega_{S}} + \tilde{A}(v) + \tilde{B}(v)\;.       
\end{equation}
Here, $A(v)$ and $\tilde{A}(v)$ consist of terms of order $k$ with at least one of the derivatives being barred or complex tangential, $B(v)$ and $\tilde{B}(v)$ consist of terms of order at most $(k-1)$. The general case follows by induction. The error terms $A(v)$ and $\tilde{A}(v)$ are benign, again by \eqref{benign1} and \eqref{benign2}. So $\alpha_{\rho}$ governs the relevant commutators in the $\overline{\partial}$--Neumann problem.

\smallskip

It is important in section \ref{alphaest}  that certain properties of $\alpha_{\rho}$ do not depend on the particular choice of $\rho$. This independence stems from the following observation. When $\rho$ is a fixed defining function, any other defining function can be written in the form $e^{h}\rho$, and at points of the boundary, $(L_{n})_{e^{h}\rho}=e^{-h}(L_{n})_{\rho}$. The associated D'Angelo forms then enjoy a simple relationship:
\begin{equation}\label{relation}
 \alpha_{e^{h}\rho}(\overline{L}) = \alpha_{\rho}(\overline{L}) + dh(\overline{L})\;
\end{equation}
(from \eqref{alpha1} or \eqref{alpha2}), and similarly for $L$ ($\alpha$, $h$ are real). In other words, $\alpha_{\rho}$ and $\alpha_{e^{h}\rho}$ differ on $T^{1,0}(b\Omega)\oplus T^{0,1}(b\Omega)$ by the differential of $h$.

\smallskip

A crucial property of $\alpha_{\rho}$, shown in \cite{BoasStraube93} (see also \cite{Straube10}, Lemma 5.14), that will be needed in section \ref{alphaest} is that its differential vanishes on the null space of the Levi form.
\begin{lemma}\label{closed}
 Let $\Omega$ be a bounded smooth pseudoconvex domain in $\mathbb{C}^{n}, P\in b\Omega$. Denote by $\mathcal{N}_{P} \subset T^{1,0}_{P}(b\Omega)$ the null space of the Levi form at $P$. If $X$ and $Y$ belong to $\mathcal{N}_{P}\oplus\overline{\mathcal{N}_{P}}$, then $(d\alpha_{\rho})_{P}(X, Y) =0$.
\end{lemma}
\noindent The proof of the lemma results from computing the differential in the following form:
\begin{equation}\label{differentialalpha}
 (d\alpha_{\rho})(X, Y) = -T_{\rho}(d\sigma_{\rho}(X, Y))(P)\;
\end{equation}
Type and skew Hermitian properties and polarization now allow one to reduce to the case $(d\alpha_{\rho})_{P}(X, \overline{X})= -T_{\rho}(d\sigma_{\rho}(X, \overline{X}))(P)$. The quantity that is being differentiated is the Levi form, which is nonnegative. Because $X\in\mathcal{N}_{P}$, the expression assumes an extremum at $P$, so the derivative has to vanish.

\smallskip

For more information on D'Angelo forms, see \cite{Straube10}, section 5.9, and \cite{Dall'AraMongodi21}.

\medskip\smallskip

\underline{\emph{2.3. The Diederich--Forn\ae ss index:}}
In 1977, Diederich and Forn\ae ss proved that if $\Omega$ is $C^{2}$, there is $\eta>0$ and a defining function $\rho_{\eta}$ such that $-(-\rho_{\eta})^{\eta}$ is plurisubharmonic on $\Omega$ (\cite{DF77b}); a simpler proof, but requiring $C^{3}$ boundary, was given by Range (\cite{Range81}). The supremum over such $\eta <1$ has come to be known as the Diederich--Forn\ae ss (or DF)--index of $\Omega$. The index can be seen as a measure of how close a domain is to having a plurisubharmonic defining function. Namely, composition of the defining function $\rho_{\eta}$ with the convex increasing function $x\rightarrow -(-x)^{\eta}$ for $x\rightarrow 0^{-}$ should be plurisubharmonic in $\Omega$ (near the boundary). The bigger $\eta<1$ is, the less `upwards bending` is needed. The primary motivation in \cite{DF77b} was to obtain bounded plurisubharmonic exhaustion  functions, the DF--index entered into regularity questions in the $\overline{\partial}$--Neumann problem through Kohn's work (\cite{Kohn99}). The index of the Diederich--Forn\ae ss worm domains is strictly less than one; the explicit value $\pi/(\pi+\text{total winding})$ has been calculated only surprisingly recently by Liu (\cite{Liu19}). The index is known to be one on domains that admit a defining function that is plurisubharmonic at the boundary (\cite{FornaessHerbig07, FornaessHerbig08}), domains that satisfy McNeal's property $(\tilde{P})$, hence those that satisfy Catlin's property($P$) (\cite{Sibony91}, proofs of Theorem 2.4 and 3.2, \cite{Harrington22}, Corollary 6.2)\footnote{Property $(P)$ is a general sufficient condition for compactness of the $\overline{\partial}$--Neumann operator introduced in \cite{Catlin82}; property $(\tilde{P})$ is more general and was introduced in \cite{McNeal02}. This class of domains includes in particular domains of finite type.}, domains that admit a family of `good' vector fields and have a `well behaved' set of infinite type points (\cite{Harrington19, AbdulHarr19}), and domains that are strictly pseudoconvex except for a simply connected complex submanifold in the boundary, and are strictly pseudoconvex in directions transverse to this submanifold (\cite{Liu19b}). Note that these are all domains where the $\overline{\partial}$--Neumann operator is known to be regular in Sobolev norms.

\section{D'Angelo forms and estimates}\label{alphaest}
With the notations from the previous sections, we say that the form $\alpha_{\rho}$ is approximately exact on the null space of the Levi form if the following holds (\cite{StraubeSucheston02}). There is a constant $C>0$ and for every $\varepsilon>0$ there exists a smooth real valued function $h_{\varepsilon}$ in a neighborhood $U_{\varepsilon}$ of $K$ such that $1/C\leq h_{\varepsilon}\leq C$ and $dh_{\varepsilon}|_{\mathcal{N}_{P}}=\alpha_{\rho}|_{\mathcal{N}_{P}} + O(\varepsilon)$, $P\in K$. Here $O(\varepsilon)$ denotes a $1$--form that satisfies $|O(\varepsilon)(X)|\leq \varepsilon C|X|$, $X\in\mathcal{N}_{P}$, where $C$ does not depend on $P$. In view of \eqref{relation}, if $\tilde{\rho}=e^{g}\rho$ is another defining function, then $\alpha_{\tilde{\rho}}$ is also approximately exact on the null space of the Levi form: the family $\{h_{\varepsilon}-g\}_{\varepsilon >0}$ will do. In particular, the assumption in the following theorem does not depend on the defining function $\rho$. The theorem is implicit in \cite{BoasStraube91, BoasStraube93}.
\begin{theorem}\label{BS91}
Let $\Omega$ be a smooth bounded pseudoconvex domain in $\mathbb{C}^{n}$, $\rho$ a defining function for $\Omega$, and $1\leq q\leq n$. Assume that $\alpha_{\rho}$ is approximately exact on the null space of the Levi form. Then the $\overline{\partial}$--Neumann operators $N_{q}$ are continuous on $W^{s}_{(0,q)}(\Omega)$ for $s\geq 0$. 
\end{theorem}
\begin{corollary}\label{plushdef}
Let $\Omega$ be a smooth bounded pseudoconvex domain in $\mathbb{C}^{n}$ that admits a defining function $\rho$ that is plurisubharmonic at points of the boundary. Then the $\overline{\partial}$--Neumann operators $N_{q}$, $1\leq q\leq n$, are continuous in $W^{s}_{(0,q)}(\Omega)$, for $s\geq 0$.
\end{corollary}
\begin{proof}
 The Hessian of $\rho$ is positive semidefinite on all vectors in $\mathbb{C}^{n}$ at points of the boundary. For $L\in \mathcal{N}_{p}$, \eqref{alternative} and the Cauchy--Schwarz inequality give $\alpha_{\rho}(\overline{L})(P)=0$. Then also $\alpha_{\rho}(L)(P)=0$ ($\alpha_{\rho}$ is real). Thus the assumptions in Theorem \ref{BS91} are satisfied with $h_{\varepsilon}\equiv 0$ for all $\varepsilon >0$.
\end{proof}
\noindent Via a different method, McNeal and Gallagher (\cite{Herbig-McNeal05}) showed that for regularity on $(0,q)$--forms for $q\geq q_{0}$, it suffcies that the
sum of the smallest $q_{0}$ eigenvalues of the Hessian of the defining function is nonnegative. It would be interesting to see how to adapt the $\alpha_{\rho}$--paradigm to also reflect the form level, and thus possibly obtain an approach to this result more in line with the current section. For now, it is worth noting that in Corollary \ref{plushdef} it suffices to assume that the complex Hessian of $\rho$ is positive semidefinite on the span of the null space of the Levi from and the complex normal.

Convex domains admit a defining function that is plurisubharmonic at the boundary: its Minkowski functional is convex, hence plurisubharmonic, near the boundary. 
\begin{corollary}\label{convex}
 Let $\Omega$ be a smooth bounded convex domain in $\mathbb{C}^{n}$. Then the $\overline{\partial}$--Neumann operators $N_{q}$, $1\leq q\leq n$, are continuous in $W^{s}_{(0,q)}(\Omega)$, for $s\geq 0$.
\end{corollary}
It is easy to see how Theorem \ref{BS91} results from Theorem \ref{vectorfield}. Set $X_{\varepsilon} = e^{h_{\varepsilon}}(L_{n})_{\rho}$, where $\{h_{\varepsilon}\}_{\varepsilon >0}$ is the family from the theorem.  Then $X_{\varepsilon}\rho = e^{h_{\varepsilon}}$, so is real and bounded uniformly and uniformly away from zero. Let $L$ be a tangential $(1,0)$--field of unit length and $P\in b\Omega$ with $L(P)\in \mathcal{N}_{p}$. The normal $(1,0)$--component of the commutator $[X_{\varepsilon}, \overline{L}]$ satisfies
\begin{multline}\label{field2}
\partial\rho([X_{\varepsilon}, \overline{L}])(P)=-e^{h_{\varepsilon}}\overline{L}h_{\varepsilon}(P)+ e^{h_{\varepsilon}}\partial\rho([(L_{n})_{\rho}, \overline{L}])(P) \\
= e^{h_{\varepsilon}}(-\overline{L}h_{\varepsilon}+ \alpha_{\rho}(\overline{L}))(P) = e^{h_{\varepsilon}}(-dh_{\varepsilon}(\overline{L}) + \alpha_{\rho}(\overline{L}))(P)\;.
\end{multline}
Because of the uniform bounds on $\{h_{\varepsilon}\}_{\varepsilon >0}$, the right hand side of \eqref{field2} is $O(\varepsilon)$. So the family $\{X_{\varepsilon}\}_{\varepsilon>0}$ satisfies the assumptions in Theorem \ref{vectorfield}.

\smallskip

It is natural to ask what happens if in Theorem \ref{BS91}, the assumption on approximate closedness of $\alpha_{\rho}$ only holds to within some fixed, but not arbitrarily good, accuracy. It is clear from the proof that one can then get estimates up to some finite level in the Sobolev scale. When dealing with estimates at level $k$ there appear commutators with $k$-th powers, leading to a factor $k$ in the quantities to be estimated (see \eqref{dbarcomm} and \eqref{dbar*comm})\footnote{The factor $k$ stems from the formula $[A,T^{k}]=\sum_{j=1}^{k}\binom{k}{j}T^{k-j}\underbrace{[...[A,T],T]...T]}_{j-fold}=kT^{k-1}[A,T] + \text{lower order}$, see e.g. \cite{DerridjTartakoff76}, Lemma 2, page 418; \cite{Straube10}, formula (3.54).}. Accordingly, one expects that the required accuracy should be proportional to $1/k$. Harrington and Liu (\cite{HarringtonLiu20}) have verified this expectation when $\alpha_{\rho}$ is approximately exact on all of $T^{1,0}(b\Omega)$ (rather than just on the null space of the Levi form) at weakly pseudoconvex points, but with added precision in terms of the spectrum of $\Box_{q}$. Denote by $\Sigma$ the set of weakly pseudoconvex boundary points, and fix a defining function $\rho$. 
\begin{theorem}[\cite{HarringtonLiu20}]\label{HL20}
Let $\Omega$ be a bounded pseudoconvex domain with smooth boundary in $\mathbb{C}^{n}$, and $1\leq q\leq n$. Assume there exists a real valued smooth function $h$ defined in a neighborhood $U$ of $\Sigma$ in $b\Omega$ so that 
\begin{equation}\label{approxexact}
 \sup_{\Sigma}|Lh-\alpha_{\rho}(L)| < \frac{1}{2kn\sqrt{2B\binom{n}{q}}}\;,
\end{equation}
for all unit $(1,0)$ vector fields $L$ in a neighborhood of the boundary and tangential at the boundary, where $k\in\mathbb{N}$ and $1/B$ is the infimum of the spectrum of $\Box_{q}$. Then the $\overline{\partial}$--Neumann operator for $(0,q)$--forms is bounded on $W^{k}_{(0,q)}(\Omega)$.
\end{theorem}
\noindent As above, whether \eqref{approxexact} holds is independent of the defining function $\rho$. If there is a sequence $\{h_{k}\}_{k\in\mathbb{N}}$ with $h_{k}$ satisfying \eqref{approxexact}, $N_{q}$ is regular at all levels. It is significant that no uniform bounds are required on $\{h_{k}\}_{k\in\mathbb{N}}$. I believe that Theorem \ref{HL20} should hold with \eqref{approxexact} only required for $L(P)\in\mathcal{N}_{P}$. If it does, the definition of approximately exact on the null space of the Levi form can be modified accordingly. For the main applications in this section, Theorems \ref{BS91} and \ref{trivialcohom}, uniform bounds are not an issue.

\medskip

Assume now that there is a smooth real submanifold $M$ in the boundary whose tangent space at every point is contained in the null space of the Levi from at the point. Then the restriction of $\alpha_{\rho}$ to $M$, $\alpha_{\rho}|_{M}$,  is closed and so represents a DeRham cohomology class, $[\alpha_{\rho}|_{M}]$, on $M$. In view of \eqref{relation}, this class is intrinsic, that is, it is independent of the defining function $\rho$. If the class is trivial (i.e.$\alpha_{\rho}$ is exact on $M$), global regularity holds:
\begin{theorem}\label{trivialcohom}
 Let $\Omega$ be a smooth bounded pseudoconvex domain in $\mathbb{C}^{n}$. Suppose there is a smooth real submanifold $M$ (with or without boundary) of $b\Omega$ that contains all the boundary points of infinite type and whose real tangent space at each point is contained in the null space of the Levi form at that point. If $[\alpha_{\rho}|_{M}] =0$, then the $\overline{\partial}$--Neumann operators $N_{q}$, $1\leq q\leq n$, are continuous in the Sobolev space $W_{(0,q)}^{s}(\Omega)$ when $s\geq 0$.
\end{theorem}
\noindent This is the main theorem in \cite{BoasStraube93}, see also \cite{Straube10}, Corollary 5.16. 
An important special case is the following:
\begin{corollary}\label{simplyconnected}
Let $\Omega$ be a smooth bounded pseudoconvex domain in $\mathbb{C}^{n}$, of finite type except for a smooth (as a submanifold with boundary) simply connected complex manifold in the boundary. Then the $\overline{\partial}$--Neumann operators $N_{q}$, $1\leq q\leq n$, are continuous in $W_{(0,q)}^{s}(\Omega)$, $s\geq 0$.
\end{corollary}
\noindent Note that such a complex submanifold is a closed, hence compact, subset of the boundary and so must be a submanifold with boundary.

I now sketch the proof of Theorem \ref{trivialcohom}, for details, see \cite{BoasStraube93} or \cite{Straube10}, Corollary 5.17. Because $[\alpha_{\rho}|_{M}]=0$, there is a smooth function $h$ on $M$ such that $d_{M}h=\alpha_{\rho}|_{M}$. Fix $\varepsilon>0$. By appropriately prescribing derivatives transverse to $M$ at $P$, $h$ is first extended locally, near a point $P\in M$, to obtain 
a (possibly complex valued) function $\tilde{h}_{P}$ that satisfies $\overline{L}\tilde{h}_{P}(P) = \alpha_{\rho}(\overline{L})(P)\;;\;L\in T^{1,0}(b\Omega)$. Then $\partial\rho([e^{\tilde{h}_{P}}(L_{n})_{\rho}, \overline{L}])(P)=e^{\tilde{h}_{P}}(\alpha_{\rho}(\overline{L})-\overline{L}\tilde{h}_{P})(P)=0$. In a small enough neighborhood of $P$, this relation holds with an error that is $O(\varepsilon)$, and we can also keep bounds on $|\tilde{h}_{P}|$ and $|\arg \tilde{h}_{P}|$ uniformly in $\varepsilon$. $M$ can be covered by finitely many of these neighborhoods. With a suitable partition of unity $\{\varphi_{j}\}_{j=1}^{s}$, we define $X_{\varepsilon} = \sum_{j=1}^{s}\varphi_{j}e^{\tilde{h}_{P_{j}}}(L_{n})_{\rho}$. Then the family $\{X_{\varepsilon}\}_{\varepsilon>0}$ satisfies the assumptions of Theorem \ref{vectorfield}; the partition of unity does not interfere with the commutators because at points of $K$, $\tilde{h}_{P_{j}}=h$, so that
$\sum_{j}(\overline{L}\varphi_{j})e^{\tilde{h}_{P_{j}}}(L_{n})_{\rho}=\big(\sum_{j}(\overline{L}\varphi_{j})\big)e^{h}(L_{n})_{\rho}=0$.

\smallskip

In view of the Barrett/Christ irregularity results (\cite{Barrett92, Christ96}), Theorem \ref{trivialcohom} cannot apply to worm domains. It is not hard to check that indeed $[\alpha_{\rho}|_{M}]\neq 0$. Following \cite{ChenShaw01},  section 6.4, a family of worm domains can be defined by
\begin{equation}\label{worm}
 \Omega_{r} = \{(z_{1}, z_{2})\in \mathbb{C}^{2}\,|\,|z_{1}+e^{i\log|z_{2}|^{2}}|^{2} < 1 - \varphi_{r}(\log|z_{2}|^{2}\} \;,
\end{equation}
where $\varphi_{r}$ satisfies certain conditions that ensure that $\Omega_{r}$ is a bounded smooth pseudoconvex Hartogs domain, strictly pseudoconvex except at the critical annulus $M=\{(0,z_{2})\,|\,e^{-r/2}\leq |z_{2}|\leq e^{r/2}\}$. Let $\rho_{r}$ be the obvious defining function. A computation (from \cite{Straube10}) gives, at points of $M$: 
\begin{equation}\label{alphacalc}
(L_{n})_{\rho_{r}}=e^{i\log|z_{2}|^{2}}\frac{\partial}{\partial z_{1}}\;\;\text{and}\;\;[(L_{n})_{\rho_{r}}, \frac{\partial}{\partial\overline{z_{2}}}] = -\frac{i}{\overline{z_{2}}}e^{i\log|z_{2}|^{2}}\frac{\partial}{\partial z_{1}}\;,
\end{equation}
so that $\alpha_{\rho_{r}}(\partial/\partial\overline{z_{2}}) = -i/\overline{z_{2}}$ (in view of \eqref{alpha2}). Since $\alpha_{\rho_{r}}$ is real, we obtain
\begin{equation}\label{alphaworm}
 \alpha_{\rho_{r}}|_{M} = \frac{i}{z_{2}}dz_{2} - \frac{i}{\overline{z_{2}}}d\overline{z_{2}}\;.
\end{equation}
Note that $(d\alpha_{\rho_{r}})|_{M} = d_{M}(\alpha_{\rho_{r}}|M) = 0$, in agreement with Lemma \ref{closed}. Integration of $\alpha_{\rho_{r}}$ around a concentric circle in $M$ results in $-4\pi$, so $[\alpha_{\rho_{r}}|_{M}]\neq 0$.

As mentioned in section \ref{history}, the appearance of the class $[\alpha_{\rho}|_{M}]$ explains why a disc in the boundary is benign for the $\overline{\partial}$--Neumann problem, while an annulus is bad on the worm domains, but is not on `nowhere wormlike' Hartogs domains (\cite{BoasStraube92}). It is interesting to note that in the situation of Theorem \ref{trivialcohom}, $[\alpha_{\rho}|_{M}]$ is the obstruction to the existence of the  vector fields in Theorem \ref{vectorfield}: the vector fields exist if and only if the class is trivial (\cite{BoasStraube93}, Remark 5 in section 4; \cite{BDD07}, Theorem 3.36). Whether it is likewise the obstruction to regularity is a very interesting question that is currently not understood (see however \cite{Barrett94}).

\smallskip

When there is a foliation by complex submanifolds in the boundary, say of co--dimension one, instead of just one complex submanifold as in Theorem \ref{trivialcohom}, the restriction of a D'Angelo form to each leaf is of course closed, but this property does not reflect the fact that each leaf is part of the foliation, and more can be said: the cohomology class of the restriction represents the leaf's infinitesimal holonomy in the foliation  (\cite{StraubeSucheston03}, Remark 3; \cite{BDD07}, Section 3.6). Moreover, whether the form is exact on the leaves of the foliation by a `global' function is equivalent to a question much studied in foliation theory, namely whether the foliation can be defined globally by a closed one--form (\cite{Straube06}, Section 3.1, \cite{Straube10}, Lemma 5.18, \cite{StraubeSucheston03}, Proposition 2). For a more detailed discussion of D'Angelo forms in connection with foliations by complex submanifolds in the boundary, the reader should consult \cite{StraubeSucheston03}, section 4 and especially section 3.6 of \cite{BDD07}.

\medskip

In Theorems \ref{BS91} and \ref{HL20}, the assumptions on $\alpha_{\rho}$ amount to a pointwise estimate, while in Theorem \ref{trivialcohom} the assumption concerns the DeRham cohomology class $[\alpha|_{M}]$. Yet a different type of condition was introduced in \cite{Straube05}, namely an $L^{2}$--type condition. If $u=\sideset{}{'}\sum_{|J|=q}u_{J}d\overline{z_{J}}$, and $K$ is a $(q-1)$--tuple, define the vector field $L^{K}_{u}=\sum_{j=1}^{n}u_{jK}(\partial/\partial z_{j})$. When $u\in dom(\overline{\partial}^{*})$, $L^{K}_{u}$ is tangential at the boundary.
\begin{theorem}[\cite{Straube05}]\label{St05}
 Let $\Omega$ be a smooth bounded pseudoconvex domain in $\mathbb{C}^{n}$, $\rho$ a defining function for $\Omega$. Let $1\leq q\leq n$. Assume there is a constant $C$ such that for all $\varepsilon>0$ there exists a defining function $\rho_{\varepsilon}=e^{h_{\varepsilon}}\rho$ and a constant $C_{\varepsilon}$ with $1/C < |h_{\varepsilon}|<C$ on $b\Omega$, and
\begin{equation}\label{L2cond}
 \sideset{}{'}\sum_{K}\int_{\Omega}|\alpha_{\rho_{\varepsilon}}(L^{K}_{u})|^{2} \leq \varepsilon\big(\|\overline{\partial}u\|^{2}+\|\overline{\partial}^{*}u\|^{2}\big)+C_{\varepsilon}\|u\|_{-1}^{2}\;;\;u\in C^{\infty}_{(0,q)}(\overline{\Omega})\cap dom(\overline{\partial}^{*})\;.
 \end{equation}
 Then the $\overline{\partial}$--Neumann operator $N_{q}$ is continuous on $W^{s}_{(0,q)}(\Omega)$, $s\geq 0$. 
 \end{theorem}
\noindent The motivation in \cite{Straube05} was to find an approach to regularity that combines the vector field method with that via compactness, and for that the uniform bounds required on $\{h_{\varepsilon}\}$ were acceptable at the time.\footnote{For another generalization of the vector field method that incorporates potential theoretic conditions, see \cite{Harrington11}.} In connection with the DF--index, such bounds are not available. But the step from working with pointwise estimates of some kind on $\alpha_{\rho_{\epsilon}}$ to an $L^{2}$--type estimate is nonetheless essential in section \ref{DF-reg}. And it is entirely conceivable that Theorem \ref{St05} holds without requiring the uniform bounds. In fact, the work in section \ref{DF-reg} (from \cite{LiuStraube22}) shows that it holds in this generality on domains with comparable Levi eigenvalues (see section \ref{DF-reg} for the definition), in particular for domains in $\mathbb{C}^{2}$.

\section{D'Angelo forms and Diederich--Forn\ae ss index}\label{DA-DF}
In this section, I discuss results by Liu(\cite{Liu19}) and Yum (\cite{Yum21}) (see also \cite{Dall'AraMongodi21}) that relate the DF--index with D'Angelo forms.

Let $\eta\in (0,1)$ and assume $\rho_{\eta}$ is a defining function for $\Omega$ such that $-(-\rho_{\eta})^{\eta}$ is plurisubharmonic in $\Omega$, close to the boundary. Computing the Hessian of $-(-\rho_{\eta})^{\eta}$ shows that for $z$ close enough to the boundary (which from now on is assumed)
\begin{equation}\label{defA}
 A_{\eta}(w, \overline{w}) := \sum_{j,k=1}^{n}\left(\frac{\partial^{2}\rho_{\eta}}{\partial z_{j}\partial\overline{z_{k}}} + \frac{1-\eta}{-\rho_{\eta}}\,\frac{\partial\rho_{\eta}}{\partial z_{j}}\,\frac{\partial\rho_{\eta}}{\partial\overline{z_{k}}}\right)(z)\,w_{j}\overline{w_{k}} \geq 0\;,\;w\in\mathbb{C}^{n}\;.
\end{equation}
Recall from \eqref{alternative} that if $L=\sum_{j=1}^{n}w_{j}(\partial/\partial z_{j})$, then 
$\alpha_{\rho_{\eta}}(\overline{L}) = \sum_{j,k=1}^{n}\frac{\partial^{2}\rho_{\eta}}{\partial z_{j}\partial\overline{z_{k}}}((L_{n})_{\rho_{\eta}})_{j}\overline{w_{k}}$, that is, it is the mixed term in the Hessian in \eqref{defA}. If $L$ is tangential to the level surfaces of $\rho_{\eta}$, this mixed term agrees with $A_{\eta}((L_{n})_{\rho_{\eta}}, \overline{L})$, and we can use Cauchy--Schwarz to estimate it:
\begin{equation}\label{CS}
|\alpha_{\eta}(\overline{L})|^{2} = |A_{\eta}((L_{n})_{\rho_{\eta}}, \overline{L})|^{2}
\leq A_{\eta}\left((L_{n})_{\rho_{\eta}}, \overline{(L_{n})_{\rho_{\eta}}}\right)A_{\eta}(L, \overline{L})\;.
\end{equation}
For $P\in b\Omega$, $L(P)\in \mathcal{N}_{P}$, the idea is now to obtain an estimate on $\alpha_{\eta}(\overline{L})(P)$ by letting $z$ tend to $P$ along the integral curve of the gradient of $\rho_{\eta}$ that goes through $P$ and computing the limit of the right hand side of \eqref{CS}\footnote{This idea also occurs in \cite{Straube01}, where it is used in connection with the plurisubharmonicity of $-\log(-\rho)$ when $(\rho)$ is the boundary distance, see estimate (6) and the subsequent argument there.}. We have 
\begin{equation}\label{ineq1}
A_{\eta}\left((L_{n})_{\rho_{\eta}}, \overline{(L_{n})_{\rho_{\eta}}}\right)(z) = \frac{1-\eta}{\rho_{\eta}}\left(\frac{\rho_{\eta}}{1-\eta}\sum_{j,k=1}^{n}\frac{\partial^{2}\rho_{\eta}}{\partial z_{j}\partial\overline{z_{k}}}((L_{N})_{\rho_{\eta}})_{j}\overline{((L_{N})_{\rho_{\eta}})_{k}} \;\;- 1\right)(z)\;.
\end{equation}
The expression in the parenthesis on the right tends to $(-1)$ as $z\rightarrow P$. Therefore we are left with the limit of $((1-\eta)/\rho_{\eta})A_{\eta}(L, \overline{L})=((1-\eta)/\rho_{\eta})\sum_{j,k=1}^{n}\frac{\partial^{2}\rho_{\eta}}{\partial z_{j}\partial\overline{z_{k}}}w_{j}\overline{w_{k}}$ ($L$ is tangent to the level surfaces of $\rho_{\eta}$). Because $L(P)\in\mathcal{N}_{P}$, this is a (real) normal derivative. Also, when $L(P)\in\mathcal{N}_{P}$, tangential derivatives vanish (as in the proof of Lemma \ref{closed}), and because we are dividing by $\rho_{\eta}$ rather than by minus the boundary distance, we end up with $(L_{n})_{\rho_{\eta}}\big(\sum_{j,k=1}^{n}\frac{\partial^{2}\rho_{\eta}}{\partial z_{j}\partial\overline{z_{k}}}w_{j}\overline{w_{k}}\big)(P)$ (compare Lemma 1.1 in \cite{Liu19}). Computing this derivative is the key step in the argument; this computation is done in \cite{Liu19}, Lemmas 1.1 and 2.1, and in \cite{Dall'AraMongodi21}, Lemma 6.8. Related expressions are computed in \cite{Yum21}, Proposition 4.6. Following  \cite{Dall'AraMongodi21}, observe that $\sum_{j,k=1}^{n}\frac{\partial^{2}\rho_{\eta}}{\partial z_{j}\partial\overline{z_{k}}}w_{j}\overline{w_{k}}=2\partial\overline{\partial}\rho_{\eta}(L, \overline{L})$, and that the form $\partial\overline{\partial}\rho_{\eta}$ is closed; in particular, this means that $d\partial\overline{\partial}\rho_{\eta}\big((L_{n})_{\rho_{\eta}}, L, \overline{L}\big)=0$. Writing out this differential with the usual formula\footnote{I use the conventions from \cite{Morita01}.} gives six terms, one of which is the desired derivative. Of the remaining five terms, two vanish for simple reasons, and the other three can all be expressed in terms  of $\alpha_{\rho_{\eta}}$. One again needs to use that $L(P)\in\mathcal{N}_{P}$ (see \cite{Dall'AraMongodi21} for details). Solving for the term we are interested in results in 
\begin{multline}\label{derivative}
 (L_{n})_{\rho_{\eta}}\big(\sum_{j,k=1}^{n}\frac{\partial^{2}\rho_{\eta}}{\partial z_{j}\partial\overline{z_{k}}}w_{j}\overline{w_{k}}\big)(P) \\
 = L(\alpha_{\rho_{\eta}}(\overline{L}))(P)-|\alpha_{\rho_{\eta}}(L)(P)|^{2}-\alpha_{\rho_{\eta}}([L, \overline{L}]_{(0,1)})(P)\;,\;L(P)\in\mathcal{N}_{P}\;.
\end{multline}
Subscripts denote the type of a form or of a vector field. Type considerations give (compare also Lemma 4.6 in \cite{Dall'AraMongodi21})
\begin{multline}\label{1}
 L(\alpha_{\rho_{\eta}}(\overline{L}))-\alpha_{\rho_{\eta}}([L, \overline{L}]_{(0,1)})\\
= L\big((\alpha_{\rho_{\eta}})_{(0,1)}(\overline{L})\big)-\overline{L}\big((\alpha_{\rho_{\eta}})_{(0,1)}(L)\big)-(\alpha_{\rho_{\eta}})_{(0,1)}([L, \overline{L}])\\
= 2d\big((\alpha_{\rho_{\eta}})_{(0,1)}\big)(L, \overline{L})) = 2\partial\big((\alpha_{\rho_{\eta}})_{(0,1)}\big)(L,\overline{L}) = 2\partial\alpha_{\rho_{\eta}}(L, \overline{L})\;.
\end{multline}
Combining \eqref{CS}, the limit $(-1)$ for the term in the parenthesis on the right hand side of \eqref{ineq1}, \eqref{derivative}, and \eqref{1}  gives
\begin{equation}\label{ineq2}
 |\alpha_{\rho_{\eta}}(L)(P)|^{2} \leq (1-\eta)\big(-2\partial\alpha_{\rho_{\eta}}(L, \overline{L})(P) + |\alpha_{\rho_{\eta}}(L)(P)|^{2}\big)\;,
\end{equation}
or, since $L(P)\in\mathcal{N}_{P}$ and Lemma \ref{closed} imply $\partial\alpha_{\rho_{\eta}}(L, \overline{L})(P)+ \overline{\partial}\alpha_{\rho_{\eta}}(L, \overline{L})(P)= 0$,
\begin{equation}\label{ineq3}
 |\alpha_{\rho_{\eta}}(L)(P)|^{2} \leq 2\frac{(1-\eta)}{\eta}\overline{\partial}\alpha_{\rho_{\eta}}(L, \overline{L})(P)\;.
\end{equation}
This inequality establishes the $\leq$ direction  in the following theorem (\cite{Liu19}, Theorem 2.10; \cite{Yum21}, Theorem 1.1):
\begin{theorem}\label{DFDA2}
$DF(\Omega) = \sup\{\eta\in (0,1)|\eqref{ineq3}\;\text{holds for all $P\in b\Omega$, $L(P)\in\mathcal{N}_{P}$}\}$. 
\end{theorem}
\noindent Note that the left hand side of \eqref{ineq3} equals $2\big((\alpha_{\rho_{\eta}})_{(0,1)}\wedge\overline{(\alpha_{\rho_{\eta}})_{(0,1)}})(L, \overline{L})$, so that \eqref{ineq3} agrees with the inequality given in \cite{Yum21}, Theorem 1.1. An equivalent version of \eqref{ineq3}, and hence of the quantity to take the supremum over in Theorem \ref{DFDA2}, is obtained by not computing the derivative in \eqref{derivative} and instead leaving it in the estimate, and expressing $\alpha_{\rho_{\eta}}(L)$ via \eqref{alternative}:
\begin{multline}\label{DFDA3}
\frac{1}{1-\eta}\big(\sum_{j,k=1}^{n}\frac{\partial^{2}\rho_{\eta}}{\partial z_{j}\partial\overline{z_{k}}}((L_{n})_{\rho_{\eta}})_{j}\overline{w_{k}}\,\,\big)^{2}(P) + (L_{n})_{\rho_{\eta}}\big(\sum_{j,k=1}^{n}\frac{\partial^{2}\rho_{\eta}}{\partial z_{j}\partial\overline{z_{k}}}w_{j}\overline{w_{k}}\,\big)(P) \leq 0\;;  \\
L(P)\in\mathcal{N}_{P}\;.
\end{multline}
\noindent Yum indicated in \cite{Yum21}, Theorem 4.1 and Remark 4.3 that \eqref{DFDA3} is essentially the inequality in \cite{Liu19} that is implied by Theorem 2.10 there. He also realized that \eqref{DFDA3} can be neatly expressed via D'Angelo forms as in \eqref{ineq3}, and Dall'Ara and Mongodi (\cite{Dall'AraMongodi21}) sorted things out nicely in coordinate free language.
\smallskip

To prove the $\geq$ direction in Theorem \ref{DFDA2}, let $\eta$ and $\rho_{\eta}$ be so that \eqref{ineq3} holds. (Because $0< DF(\Omega)$ and the direction already proved, the set of these $\eta$'s is not empty.) It suffices to show that for any $\delta$ with $0<\delta<\eta$, there is a defining function $\rho_{\delta}$ such that $-(-\rho_{\delta})^{\delta}$ is plurisubharmonic in $\Omega$ (near the boundary\footnote{ See Lemma 2.3 in \cite{Harrington22}.}), equivalently $A_{\delta}(w, \overline{w}) \geq 0$, $w\in\mathbb{C}^{n}$. The argument actually establishes strict definiteness. It suffices to show that the Hermitian form in $a, b\in\mathbb{C}$ given by $A_{\rho_{\delta}}(aL+b(L_{n})_{\rho_{\delta}}, \overline{aL+b(L_{n})_{\rho_{\delta}}})$ is positive definite whenever $L$ is tangent to the level sets of $\rho_{\delta}$. The rest of this proof uses ideas from \cite{Dall'AraMongodi21}, proof of Theorem 6.6; similar computations can also be found in \cite{Harrington22}. Looking at the minors of the determinant of this Hermitian form in order to apply Sylvester's criterion, note that $A_{\rho_{\delta}}\big((L_{n})_{\rho_{\delta}}, \overline{(L_{n})_{\rho_{\delta}}})\geq (1/2)(1-\delta)/(-\rho_{\delta})>0$ close enough to the boundary. Computation of the determinant and multiplication by $(-\rho_{\delta})$ reveals that the condition for its positivity is
\begin{multline}\label{detcomp}
\;\;\;\;\; (-\rho_{\delta})\left((2\partial\overline{\partial}\rho_{\delta})(L,\overline{L})(2\partial\overline{\partial}\rho_{\delta})((L_{n})_{\rho_{\delta}}, \overline{(L_{n})_{\rho_{\delta}}}) - |(2\partial\overline{\partial}\rho_{\delta})(L,\overline{(L_{n})_{\rho_{\delta}}})|^{2}\right) \\
 +(1-\delta)(2\partial\overline{\partial}\rho_{\delta})(L,\overline{L}) > 0\;;\;\;\;\;\;
\end{multline}
recall that $(2\partial\overline{\partial}\rho_{\delta})$ corresponds to the Hessian of $\rho_{\delta}$. On the boundary, the left hand side of \eqref{detcomp} equals $(1-\delta)(2\partial\overline{\partial}\rho_{\delta})(L,\overline{L})\geq 0$.  Regardless of how $\rho_{\delta}$ is chosen, if at a point $\partial\overline{\partial}\rho_{\delta}(L,\overline{L})$ is strictly positive, then the left hand side of \eqref{detcomp} is strictly positive in a neighborhood of that point (inside $\Omega$). So it suffices to find $\rho_{\delta}$ so that this latter property also holds at points where $\partial\overline{\partial}\rho_{\delta}(L,\overline{L})$ vanishes. It is not hard to see that this will happen if the inside normal derivative of the left hand side of \eqref{detcomp} is strictly positive\footnote{Only requiring nonnegativity does not work. So if this argument is to establish plurisubharmonicity of $-(-\rho_{\eta})^{\eta}$, it actually must establish strict plurisubharmonicity.}, but see \cite{Dall'AraMongodi21}, Lemma 6.9. Because the unit sphere bundle in $T^{1,0}(b\Omega)$ is compact, we then find a neighborhood that works for all sections.

Set $\rho_{\delta}:=e^{-\varepsilon|z|^{2}}\rho_{\eta}$, where $\varepsilon$ will be chosen sufficiently small later, and $\rho_{\eta}$ is from above. Let $P\in b\Omega$ with $\partial\overline{\partial}\rho_{\delta}(L,\overline{L})(P)=0$. The derivative of the left hand side of \eqref{detcomp} at $P$ equals
\begin{equation}\label{valentine1}
 -\delta|\alpha_{\rho_{\delta}}(L)|^{2}+2(1-\delta)\overline{\partial}\alpha_{\rho_{\delta}}(L,\overline{L})\;;
\end{equation}
I have used \eqref{alternative}, and \eqref{derivative} and \eqref{1} again for the derivative of $(2\partial\overline{\partial}\rho_{\delta})(L,\overline{L})$. We want \eqref{valentine1} to be strictly positive at $P$. Now $\alpha_{\rho_{\delta}}(L)=\alpha_{\rho_{\eta}}(L)-\varepsilon d(|z|^{2})(L)$ (by \eqref{relation}) and $\overline{\partial}\alpha_{\rho_{\delta}}(L,\overline{L})=\overline{\partial}\alpha_{\rho_{\eta}}(L,\overline{L})+\varepsilon \partial\overline{\partial}(|z|^{2})(L,\overline{L})$ (see the computation right after Corollary \ref{n=2} below). Therefore,
\begin{equation}\label{valentine3}
 -\delta|\alpha_{\rho_{\delta}}(L)|^{2}=-\delta\left(|\alpha_{\rho_{\eta}}(L)|^{2}-2\varepsilon Re\alpha_{\rho_{\eta}}(L)(d|z|^{2})(L,\overline{L})+\varepsilon^{2}|d(|z|^{2})(L)|^{2}\right) \;,
 \end{equation}
 and
 \begin{equation}\label{valentine4}
2(1-\delta)\overline{\partial}\alpha_{\rho_{\delta}}(L,\overline{L})=2(1-\delta)\big(\overline{\partial}\alpha_{\rho_{\eta}}(L,\overline{L})+\varepsilon \partial\overline{\partial}(|z|^{2})(L,\overline{L})\big)\;.
 \end{equation}
Adding the right hand sides of \eqref{valentine3} and \eqref{valentine4} it is easy to see that one can choose $\varepsilon>0$ so that \eqref{valentine1} is strictly positive at $P$.

It is noteworthy that in the construction above, $-(-\rho_{\delta})^{\delta}$ is actually strictly plurisubharmonic near the boundary. In particular, one may restrict attention to strictly plurisubharmonic functions in the definition of the DF--index (compare also \cite{Harrington22}).

\medskip

In view of \eqref{alpha2}, \eqref{alpha3}, \eqref{dbarcomm}, and \eqref{dbar*comm}, D'Angelo forms link the Diederich--Fornaess index to commutators via Theorem \ref{DFDA2} (specifically through the estimate \eqref{ineq3}), and therefore to regularity estimates, at least in principle.  In the next section, I discuss situations where this connection can indeed be exploited to obtain estimates.

\section{DF--index one and regularity}\label{DF-reg}

As of this writing, it is open whether or not DF--index one implies regularity in Sobolev spaces for the $\overline{\partial}$--Neumann operators. So far, conditions are needed. The results by Kohn (\cite{Kohn99}), then improved by Harrington (\cite{Harrington11}), impose a condition on the growth of $|\alpha_{\rho_{\eta}}|$ as $\eta\rightarrow 1^{-}$, where $\rho_{\eta}$ is the defining function from the DF--index. As always, the action is on the null space of the Levi from, and that is where the growth is taken. Recall that $\Sigma$ denotes the set of weakly pseudoconvex boundary points and define $m_{\eta}=\sup_{P\in\Sigma}\sup_{X\in\mathcal{N}_{P}\,;\,|X|=1}|\alpha_{\rho_{\eta}}(X)|$. The result is the following.
\begin{theorem}[\cite{Harrington11}]\label{Har11}
Let $\Omega$ be a smooth bounded pseudoconvex domain in $\mathbb{C}^{n}$. Suppose that for every $\eta$ in $(0,1)$, there exists a defining function $\rho_{\eta}$ such that $-(-\rho_{\eta})^{\eta}$ is plurisubharmonic and 
\begin{equation}\label{liminf}
 \liminf_{\eta\rightarrow 1^{-}}\sqrt{1-\eta}\,(m_{\eta}) = 0\;.
\end{equation}
Then the $\overline{\partial}$--Neumann operators $N_{q}$, $1\leq q\leq n$, are continuous in $W^{s}_{(0,q)}(\Omega)$, $s\geq 0$.
\end{theorem}
\noindent Kohn's result in \cite{Kohn99} had the factor $(1-\eta)^{1/3}$ in \eqref{liminf}.
Theorem \ref{Har11} is a consequence of Theorem 5.1 and Theorem 6.2 in \cite{Harrington11}\footnote{The formulation in \cite{Harrington11} uses $|Xh_{\eta}|$, where $\rho_{\eta}=e^{h_{\eta}}\rho$, instead of $|\alpha_{\rho_{\eta}}(X)|$ in the definition of $m_{\eta}$. Because $\alpha_{\rho_{\eta}}=\alpha_{\rho}+dh_{\eta}$ on $T^{1,0}(b\Omega)$, the resulting conditions from \eqref{liminf} are equivalent.}. Harrington proves that \eqref{liminf} implies that the assumptions in his generalization of the vector field method in \cite{Harrington11} (alluded to in footnote 12 above) are satisfied.

The first result that does not impose restrictions on the growth of the D'Angelo forms $\alpha_{\rho_{\eta}}$ associated with the defining functions $\rho_{\eta}$ in the definition of the index was found by Liu in \cite{Liu19b}. Instead, he assumed that the Levi eigenvalues are comparable. In addition, he needed an additional ad hoc and somewhat restrictive geometric condition. In \cite{LiuStraube22}, we were able to dispense with this latter condition.

The condition on the eigenvalues of the Levi form is the following: sums of any $q$ eigenvalues ($q$--sums) of the Levi form should dominate the trace of the Levi form. Because all eigenvalues are nonnegative, it is easy to see that this condition is equivalent to saying that any two $q$--sums should be comparable (uniformly over the boundary). Also, if the condition holds at level $q$, it holds at level $(q+1)$. The condition is independent of the choice of basis, even if orthogonality is dropped (thanks to a theorem of Ostrowski, \cite{HornJohnson85}, Theorem 4.5.9). We will need the assumption on the Levi eigenvalues in section \ref{DF-reg} because it implies (in fact, is equivalent to) so called maximal estimates (\cite{Derridj78}, Th\'{e}or\`{e}me 3.1 for $q=1$; \cite{BenMoussa00}, Th\'{e}or\`{e}me 3.7 for $q>1$).  Let $Y$ be $(1,0)$ vector field, smooth on $\overline{\Omega}$, and tangential at the boundary. A maximal estimate on $(0,q)$--forms means the estimate
\begin{equation}\label{maxest}
 \|Yu\|_{k-1}^{2} \leq C\big(\|\overline{\partial}u\|_{k-1}^{2} + \|\overline{\partial}^{*}u\|_{k-1}^{2} + \|u\|_{k-1}^{2}\big)\;;\;u\in C^{\infty}_{(0,q)}(\overline{\Omega})\;.
\end{equation}
This estimate is slightly, but decisively, stronger than \eqref{benign2}, where $\|u\|_{k-1}^{2}$ is replaced by $\|u\|_{k-1}\,\|u\|_{k}$. For more information on maximal estimates and their role in the $\overline{\partial}$--Neumann theory, see the introduction in \cite{Koenig15}, compare also \cite{Derridj91}. It is worthwhile to note that the condition of comparable sums of Levi eigenvalues is purely in terms of the boundary geometry of $\Omega$.

Then we have
\begin{theorem}[\cite{LiuStraube22}]\label{LiuStr}
 Let $\Omega$ be a smooth bounded pseudoconvex domain in $\mathbb{C}^{n}$, $1\leq q_{0}\leq (n-1)$. Assume that $q_{0}$--sums of the eigenvalues of the Levi form are comparable. Then if the Diederich--Forn\ae ss index of $\Omega$ is one, the $\overline{\partial}$--Neumann operators $N_{q}$, $q_{0}\leq q\leq n$, are continuous in Sobolev $s$--norms for $s\geq 0$.
\end{theorem}
\noindent When $q=(n-1)$, the comparability condition is trivially satisfied, as there is only one $(n-1)$--sum.
\begin{corollary}\label{n-1}
 Let $\Omega$ be a smooth bounded pseudoconvex domain in $\mathbb{C}^{n}$ with Diederich--Forn\ae ss index one. Then $N_{n-1}$ is continuous in Sobolev $s$--norms for $s\geq 0$.
\end{corollary}
\noindent I single out the important case $n=2$:
\begin{corollary}\label{n=2}
 Let $\Omega$ be a smooth bounded pseudoconvex domain in $\mathbb{C}^{2}$ with Diederich--Forn\ae ss index one. Then $N_{1}$ is continuous in Sobolev $s$--norms for $s\geq 0$.
\end{corollary}

The proof of Theorem \ref{LiuStr} results from the following ideas. Assume $\eta$ and $\rho_{\eta}$ are such that \eqref{ineq3} holds. Fix some defining function $\rho$, and define $h_{\eta}$ (near the boundary) via $\rho_{\eta}=e^{h_{\eta}}\rho$. Let $L(P)\in\mathcal{N}_{P}$. In view of Lemma \ref{closed} and of \eqref{1}, we have
\begin{multline}\label{spring1}
2\overline{\partial}\alpha_{\rho_{\eta}}(L,\overline{L})(P)=-2\partial\alpha_{\rho_{\eta}}(L,\overline{L})(P)=-L\big(\alpha_{\rho_{\eta}}(\overline{L})\big)(P)+\alpha_{\rho_{\eta}}\big([L,\overline{L}]_{(0,1)}\big)(P) \\
= -L(\alpha_{\rho}(\overline{L})+\overline{L}h_{\eta})+\alpha_{\rho}([L,\overline{L}]_{(0,1)})(P)+\big([L,\overline{L}]_{(0,1)}h_{\eta}\big)(P) \\
= 2\overline{\partial}\alpha_{\rho}(L,\overline{L})(P)-L\overline{L}h_{\eta}(P)+\big([L,\overline{L}]_{(0,1)}h_{\eta}\big)(P)\;.
\end{multline}
To pass from the first line in \eqref{spring1} to the second, use \eqref{relation} and the fact that because $L(P)\in\mathcal{N}_{(P)}$, $[L,\overline{L}](P)\in T_{P}^{1,0}\oplus T_{P}^{0,1}$. A computation shows that $-L\overline{L}h_{\eta}(P)+\big([L,\overline{L}]_{(0,1)}h_{\eta}\big)(P)=\left(\sum_{j,k}\frac{\partial^{2}(-h_{\eta})}{\partial z_{j}\partial\overline{z_{k}}}L_{j}\overline{L_{k}}\right)(P)$. Now we are in a position to use \eqref{ineq3} to obtain, for $L(P)\in\mathcal{N}_{P}$:
\begin{equation}\label{valentine7}
 |\alpha_{\rho_{\eta}}(L)(P)|^{2} \leq C(1-\eta)\big(\sum_{j,k=1}^{n}\frac{\partial^{2}(-h_{\eta})}{\partial z_{j}\partial\overline{z_{k}}}L_{j}\overline{L_{k}} + |L|^{2}\big)\;\,;\;L(P)\in\mathcal{N}_{P}\;.
\end{equation}
The $|L|^{2}$ term is needed because of the $\overline{\partial}\alpha_{\rho}$ term. Now the idea is that if $L\notin\mathcal{N}_{P}$, then $|L|$, hence $|\alpha_{\rho_{\eta}}(L)|^{2}$, can be estimated by the Levi form. This takes some care; one first extends \eqref{valentine7} into a neighborhood of the unit sphere bundle in $T^{1,0}(b\Omega)\cap\mathcal{N}$, and from there into a neighborhood $U_{\eta}$ of he boundary, see \cite{LiuStraube22}, inequality (12). This adds a term  $M_{\eta}\sum_{j,k=1}^{n}\frac{\partial^{2}\rho}{\partial z_{j}\partial\overline{z_{k}}}L_{j}\overline{l_{k}}$ to the right hand side of \eqref{valentine7}. Integrating both sides of the resulting inequality over $\Omega$, controlling compactly supported terms by interior elliptic regularity, using that integrating the Levi form over $\Omega$ rather than just over the boundary acts like a subelliptic multiplier (see \cite{Kohn79}, Proposition 4.7, part (C); \cite{D'Angelo93}, Proposition 4 in section 6.4.2), and summing over $K$, one obtains
\begin{multline}\label{valentine9}
 \sideset{}{'}\sum_{|K|=q-1}\int_{\Omega}|\alpha_{\rho_{\eta}}(L^{K}_{u})|^{2} \lesssim (1-\eta)\left(\sideset{}{'}\sum_{|K|=q-1}\int_{\Omega}\sum_{j,k=1}^{n}\frac{\partial^{2}(-h_{\eta})}{\partial z_{j}\partial\overline{z_{k}}}u_{jK}\overline{u_{kK}}+\|\overline{\partial}u\|^{2}+\|\overline{\partial}^{*}u\|^{2}\right) \\\;\;\;\;+ C_{\eta}\|u\|_{-1}^{2}\;;\;u\in C^{\infty}_{(0,q)}(\overline{\Omega})\cap dom(\overline{\partial}^{*})\;.
\end{multline}
Here the $L^{K}_{u}$ are the vector fields introduced at the end of section \ref{alphaest}, that is $L^{K}_{u}=\sum_{j=1}^{n}u_{jK}(\partial/\partial z_{j})$.
For a full explanation, see the derivation of inequality (13) in \cite{LiuStraube22}. Inequality \eqref{valentine9} holds in general, that is, without assuming comparable Levi eigenvalue sums. Next, we prove the analogue of estimate \eqref{L2cond} in Theorem \ref{St05} (without the uniform bounds on $h_{\eta}$), see Proposition 1 in \cite{LiuStraube22}:
\begin{multline}\label{prop1}
 \sideset{}{'}\sum_{|K|=q-1}\int_{\Omega}|\alpha_{\rho_{\eta}}(L^{K}_{u})|^{2} \leq C(1-\eta)\big(\|\overline{\partial}u\|^{2}+\|\overline{\partial}^{*}u\|^{2}\big) + C_{\eta}\|u\|_{-1}^{2}\;; \\
 u \in C^{\infty}_{(0,q)}(\overline{\Omega})\cap dom(\overline{\partial}^{*})\;.
\end{multline}
This estimate is obtained form \eqref{valentine9} by using the Kohn--Morrey--H\"{o}rmander formula on the right hand side, with weight $e^{h_{\eta}}$, but applied not to $u$, but to  $e^{-h_{\eta}/2}u$. Ignoring technicalities, this leads to terms $\|\overline{\partial}(e^{-h_{\eta}/2}u)\|_{-h_{\eta}}^{2}$ and $\|\overline{\partial}^{*}(e^{-h_{\eta}/2}u)\|_{-h_{\eta}}^{2}$ ($\|\cdot\|_{-h_{\eta}}$ as usual denotes norms with weight $e^{h_{\eta}}$). The latter term can be handled by essentially a self bounded gradient argument for $h_{\eta}$ (in view of the extended version of \eqref{valentine7} and $\alpha_{\rho_{\eta}}=\alpha_{\rho}+dh_{\eta}$ on $T^{1,0}(b\Omega)$; compare \cite{Straube10}, proof of Theorem 4.29 for the self bounded gradient argument). The troubling term is $\|\overline{\partial}(e^{-h_{\eta}/2}u)\|_{-h_{\eta}}^{2}$; it leads to a term $\|\overline{\partial}h_{\eta}\wedge u\|^{2}$ that we do not know how to handle without maximal estimates.

Now that we have \eqref{prop1}, the ghist of the argument in the rest of the proof of Theorem \ref{LiuStr} is as follows. Say we want to prove estimates in $W^{1}_{(0,q)}(\Omega)$. We start with showing an (a priori) estimate for $\|N_{q}u\|_{1}^{2}$. In view of \eqref{benign1} and \eqref{benign2}, we only have to estimate $\|T_{\rho_{\eta}}N_{q}u\|^{2}$. We adopt the usual conventions, in particular, $T_{\rho_{\eta}}$ acts in special boundary charts so that it preserves $dom(\overline{\partial}^{*})$. This introduces an error of order zero, which in turn contributes an error term in the estimates that can be controlled. Because the error to $T_{\rho_{\eta}}$ is of order zero, the same is true for $[\overline{\partial},T_{\rho_{\eta}}]$ and $[\overline{\partial}^{*},T_{\rho_{\eta}}]$. Terms that one can control `as usual' are denoted by `ok', and $\lesssim$ denotes an estimate with constants independent of $\eta$. 
Then, in view of \eqref{basic},
\begin{multline}\label{apriori0}
\|T_{\rho_{\eta}}N_{q}u\|^{2}\lesssim \|\overline{\partial}T_{\rho_{\eta}}N_{q}u\|^{2}+\|\overline{\partial}^{*}T_{\rho_{\eta}}N_{q}u\|^{2} \\
\lesssim \|T_{\rho_{\eta}}\overline{\partial}N_{q}u\|^{2}+\|T_{\rho_{\eta}}\overline{\partial}^{*}N_{q}u\|^{2}+\|[\overline{\partial},T_{\rho_{\eta}}]N_{q}u\|^{2}+\|[\overline{\partial}^{*},T_{\rho_{\eta}}]N_{q}u\|^{2}\;.
\end{multline}
On the first two terms on the right hand side, one uses standard integration by parts and commutations to make $\Box_{q}N_{q}u=u$ appear. Then s.c.--l.c. estimates, and absorption of terms lead to
\begin{equation}\label{apriori4}
 \|\overline{\partial}T_{\rho_{\eta}}N_{q}u\|^{2}+\|\overline{\partial}^{*}T_{\rho_{\eta}}N_{q}u\|^{2} 
 \lesssim \|T_{\rho_{\eta}}u\|^{2}+\|[\overline{\partial},T_{\rho_{\eta}}]N_{q}u\|^{2}+\|[\overline{\partial}^{*},T_{\rho_{\eta}}]N_{q}u\|^{2}+\;ok\;.
\end{equation}
By \eqref{dbarcomm}, a typical term in $[\overline{\partial},T_{\rho_{\eta}}]N_{q}u$ in a special boundary chart is $\alpha_{\rho_{\eta}}(\overline{L_{j}})(T_{\rho_{\eta}}N_{q}u)_{J}(\overline{\omega}_{j}\wedge\overline{\omega_{J}})$. So what needs to be estimated is $\int_{\Omega}|\alpha_{\rho_{\eta}}(\overline{L_{j}})T_{\rho_{\eta}}(N_{q}u)_{J}|^{2}=\int_{\Omega}|\alpha_{\rho_{\eta}}\big((T_{\rho_{\eta}}(N_{q}u)_{J})L\big)|^{2}$ ($\alpha_{\rho_{\eta}}$ is real). It is now important that we may restrict attention to $1\leq j<n$, and $n\notin J$. When $n\in J$, $T_{\rho_{\eta}}(N_{q}u)_{J}$ is taken care of because taking the normal component of a form acts like a subelliptic multiplier (\cite{Straube10}, section 2.9, \cite{Kohn79}, Proposition 4.7, part (G)). The `$j=n$`--term in \eqref{dbarcomm} can also be `neutralized' via a suitable procedure, see \cite{LiuStraube22} (reflecting the fact that commutators with the complex normal can be controlled for free, as in \cite{BoasStraube91}). So $\overline{\omega_{j}}\wedge\overline{\omega_{J}} \in dom(\overline{\partial}^{*})$ and one also checks that $\big(T_{\rho_{\eta}}(N_{q}u)_{J}\big)L_{j} = L^{J}_{(T_{\rho_{\eta}}(N_{q}u)_{J})(\overline{\omega_{j}}\wedge\overline{\omega_{J}})}$ (again because $j<n$ and $n\notin J$). Since maximal estimates (or comparable Levi sums) percolate up the form levels, they also hold for $(q+1)$--forms, and \eqref{prop1} applies (its validity does not depend on whether the vector fields $L^{J}_{u}$ are computed in Euclidean coordinates or in a special boundary chart):
\begin{multline}\label{apriori5}
 \|[\overline{\partial},T_{\rho_{\eta}}]N_{q}u\|^{2} \lesssim \sideset{}{'}\sum_{n\notin J, j<n}\int_{\Omega}|\alpha_{\rho_{\eta}}\big((T_{\rho_{\eta}}(N_{q}u)_{J})L_{j}\big)|^{2} + \;ok  \\
 \lesssim (1-\eta)\sideset{}{'}\sum_{n\notin J, j<n}\big(\|\overline{\partial}(T_{\rho_{\eta}}(N_{q}u)_{J}(\overline{\omega_{j}}\wedge\overline{\omega_{J}}))\|^{2}+\|\overline{\partial}^{*}(T_{\rho_{\eta}}(N_{q}u)_{J}(\overline{\omega_{j}}\wedge\overline{\omega_{J}}))\|^{2}\big) +\;ok\;,
\end{multline}
where $\overline{\omega_{j}}\wedge\overline{\omega_{J}}$ is suitably extended to a globally defined form (the terms on the right hand side of \eqref{apriori5} are compactly supported in the special boundary chart), and the left hand side indicates only the contribution form the particular boundary chart. The barred derivatives in the $\overline{\partial}$--term are estimated from above via \eqref{benign1} by $\|\overline{\partial}T_{\rho_{\eta}}N_{q}u\|^{2}+\|\overline{\partial}^{*}T_{\rho_{\eta}}N_{q}u\|^{2}+\|T_{\rho_{\eta}}N_{q}u\|^{2}$. For the unbarred tangential ($j<n$, $n\notin J$, so $L_{n}$ does not appear in $\overline{\partial}^{*}$) derivatives in the $\overline{\partial}^{*}$--term, \eqref{benign2} would give the term $\|T_{\rho_{\eta}}N_{q}u\|\,\|T_{\rho_{\eta}}N_{q}u\|_{1}^{2}$. The $\|\cdot\|_{1}$ term involves two derivatives of $N_{q}u$ and so cannot be absorbed (we are working on $\|N_{q}u\|_{1}$). We need the maximal estimates \eqref{maxest} once more to instead only get $\|T_{\rho_{\eta}}N_{q}u\|^{2}$. 

The commutator with $\overline{\partial}^{*}$ on the right hand side of \eqref{apriori4} is more straightforward. In \eqref{dbar*comm}, performing the summation over $j$ on the right hand shows that what must be estimated is $\sideset{}{'}\sum_{|S|=q-1}\|\alpha_{\rho_{\eta}}\big(L^{S}_{T_{\rho_{\eta}}N_{q}u}\big)\|^{2}$. Applying \eqref{prop1} once more, adding the resulting estimate to the one for the commutator with $\overline{\partial}$, and summing over boundary charts gives
\begin{multline}\label{apriori6}
 \|\overline{\partial}T_{\rho_{\eta}}N_{q}u\|^{2}+\|\overline{\partial}^{*}T_{\rho_{\eta}}N_{q}u\|^{2} \\
 \lesssim (1-\eta)\left(\|\overline{\partial}T_{\rho_{\eta}}N_{q}u\|^{2}+\|\overline{\partial}^{*}T_{\rho_{\eta}}N_{q}u\|^{2}+\|T_{\rho_{\eta}}N_{q}u\|^{2}\right)+C_{\eta}\|u\|_{1}^{2}
\end{multline}
(after also taking care of the $ok$ terms). But $\|T_{\rho_{\eta}}N_{q}u\|^{2} \lesssim \|\overline{\partial}T_{\rho_{\eta}}N_{q}u\|^{2}+\|\overline{\partial}^{*}T_{\rho_{\eta}}N_{q}u\|^{2}$ (the first inequality in \eqref{apriori0}), so choosing $\eta$ close enough to one, absorbing terms, and invoking \eqref{apriori0} again gives the desired estimate $\|T_{\rho_{\eta}}N_{q}u\|^{2}\leq C_{\eta}\|u\|_{1}^{2}$. 

\smallskip
                                                                                                                                                                                                                                                    
In \cite{LiuStraube22}, we use elliptic regularization to turn a prior estimates into genuine estimates. Doing so introduces additional technicalities into the argument above; in particular, one has to take into account that because the free boundary condition induced by the regularized quadratic form has changed, $\overline{\partial}N_{\delta,q}$ is not necessarily in the domain of $\overline{\partial}^{*}$.  It was shown in \cite{Straube05} how to handle these complications.

\medskip

The above (sketch of a) proof of Theorem \ref{LiuStr} uses DF--index one only to derive \eqref{prop1}. Afterwards, only the assumption on the Levi eigenvalues is used to derive regularity. So Theorem \ref{St05} holds without the uniform bounds on $\{h_{\varepsilon}\}$, if one assumes comparable $q$--sums of the Levi eigenvalues. In particular, the theorem holds without the bounds for domains in $\mathbb{C}^{2}$.

\vskip .5cm

\providecommand{\bysame}{\leavevmode\hbox to3em{\hrulefill}\thinspace}

\end{document}